\documentclass[journal]{IEEEtran}

\usepackage[utf8]{inputenc}
\usepackage[english]{babel}
\usepackage{amsmath}

\usepackage{multirow}
\usepackage{multicol}
\usepackage{lipsum}
\usepackage{graphicx}
\usepackage{graphics}
\usepackage{amsmath}
\usepackage{amsbsy}
\usepackage{blindtext}
\usepackage{tcolorbox}
\usepackage{amssymb}
\usepackage{scalerel}
\usepackage{mathabx}
\usepackage{verbatim}
\usepackage{booktabs}
\usepackage{rotating,tabularx}
\usepackage{mathrsfs} 
\usepackage{url}
\usepackage{adjustbox}
\usepackage{soul}
\usepackage{xcolor}
\usepackage{cite}
\usepackage{tikz}
\usepackage{color, colortbl}
\usepackage[first=0,last=9]{lcg}
\definecolor{LightGray}{rgb}{0.7,0.7,0.7}

\usepackage{algorithm}
\usepackage{algorithmic}

\usepackage{amsthm}
\theoremstyle{definition}

\theoremstyle{remark}

\usepackage[utf8]{inputenc}
\usepackage[english]{babel}

\usepackage[caption=false,font=footnotesize]{subfig}
\usepackage[T1]{fontenc}
\usepackage{scalerel,stackengine}
\newcommand\reallywidecheck[1]{%
\savestack{\tmpbox}{\stretchto{%
  \scaleto{%
    \scalerel*[\widthof{\ensuremath{#1}}]{\kern-.6pt\bigwedge\kern-.6pt}%
    {\rule[-\textheight/2]{1ex}{\textheight}}
  }{\textheight}%
}{0.5ex}}%
\stackon[1pt]{#1}{\scalebox{-1}{\tmpbox}}%
}
\stackMath
\IEEEoverridecommandlockouts
\IEEEoverridecommandlockouts

\newif\ifarxiv
 \arxivtrue

\title{Impact of Higher-Order Structures in Power Grids' Graph on Line Outage Distribution Factor}

\ifarxiv
\author{Nafis Sadik$^{\ast}$, and Mohammad Rasoul Narimani$^{\dagger}$%
\thanks{This work was supported by NSF under Award Number 2308498.}
\thanks{${\ast}$: College of Engineering and Computer Science, Arkansas State University. nafis.sadik@smail.astate.edu.}%
\thanks{$^{\dagger}$: Department of Electrical and Computer Engineering, California State University Northridge (CSUN). Rasoul.narimani@csun.edu.}}%

\else
\author{Mohammad Rasoul Narimani and Nafis Sadik}%
\fi

\begin{document}
\maketitle

\begin{abstract}
Power systems often include a specific set of lines that are crucial for the regular operations of the grid. Identifying the reasons behind the criticality of these lines is an important challenge in power system studies. When a line fails, the line outage distribution factor (LODF) quantifies the changes in power flow on the remaining lines. This paper proposes a network analysis from a local structural perspective to investigate the impact of local structural patterns in the underlying graph of power systems on the LODF of individual lines. In particular, we focus on graphlet analysis to determine the local structural properties of each line. This research analyzes potential connections between specific graphlets and the most critical lines based on their LODF. In this regard, we investigate $N-1$ and $N-2$ contingency analysis for various test cases and identifies the lines that have the greatest impact on the LODFs of other lines. We then determine which subgraphs contain the most significant lines. Our findings reveal that the most critical lines often belong to subgraphs with a less meshed but more radial structure. These findings are further validated through various test cases. Particularly, it is observed that networks with a higher percentage of ring or meshed subgraphs on their most important line (based on LODF) experience a lower LODF when that critical line is subject to an outage. Additionally, we investigate how the LODF of the most critical line varies among different test cases and examine the subgraph characteristics of those critical lines.
\end{abstract}

\begin{IEEEkeywords}
DC load flow, Line outage distribution factor, Graph theory, Graphlets. 
\end{IEEEkeywords}

\IEEEpeerreviewmaketitle

\section{Introduction}

\IEEEPARstart{T}{he} power system is comprised of interconnected networks, and the loss of even one line can potentially lead to a major blackout. As a result, it is crucial to identify lines that are critical to the network and have the potential to trigger a cascading failure and even a blackout. For that reason, remedial actions such as network redesign or transmission line capacity expansion are required\cite{galiana1984bound}. Power system contingency analysis has lately resurfaced as the traditional grid is transitioning into a smart grid. Particularly as the grid is transitioning into a cyber-physical system, grid collapse can be caused by the collapse of a few essential lines\cite{liang20162015, Boyaci_propagation}.

Numerous efforts have been expended on identifying critical lines in power systems\cite{chen2007identification, coelho2018new}. Complex network theory has been leveraged by various studies to identify critical elements in power systems. Betweenness centrality, a measure based on the shortest path-based between different nodes in a graph, is used to identify the most vulnerable lines in power system \cite{chen2007identification}. Other graph theoretical methods, such as the clustering algorithm, have also been employed to identify influential buses on the network\cite{gao2019identification}. Those research approaches focuses on the topology of the network without considering the physics that govern the network, i.e. power flow equations, line outage distribution factor, etc.
Line outage distribution factor(LODF) is an indicator that quantifies the redistribution of power flow amongst surrounding lines when a line outage happens~\cite{tejada2017security}. It is a linear sensitivity factor that utilizes DC power flow~\cite{tejada2017security}. LODF plays a vital role in optimal power flow (OPF) studies~\cite{narimani2020tightening, narimani2018comparison, narimani2018improving, narimani2018empirical, narimani2020strengthening} by quantifying the impact of line outages on power flow distribution in the network. It assists in improving system reliability, congestion management, and optimal transmission planning by providing valuable insights into the network's behavior under various operating conditions.

A combination of the line outage distribution factor and complex network analysis typically yields a better contingency analysis. To extract data from various levels of contingency, the group betweenness centrality concept can be used with the line outage distribution factor~\cite{narimani2021generalized, huang2021toward}. Nevertheless, to better analyze contingency analysis, LODFs among lines can be analyzed with the local structural properties of those lines. The literature suggests that there have been very few works that relate subgraph properties with power grid contingency analysis.

Subgraph analysis is an approach that highlights the local structural properties of a network. Graphlets are induced connected subgraph. A subgraph pattern is called a motif when it is found too frequently in a network compared to random networks\cite{milo2002network}. Motif analysis is also leveraged in power system contingency analysis studies to investigate the impact of the local structure of a network on power system characteristics. For instance, it has been discovered that robust and fragile electrical networks exhibit varying degrees of deterioration of motifs under attack\cite{dey2017motif}. However, network motifs represents more global properties rather than local properties of a network. In particular, it is a measure of how many times a subgraph is found on the whole network. This paper leveraged the local properties of power system's graph, in particular how many times a graphlet can be found in a line, in analyzing the LODF measure in power system.

The LODF is used in this paper to determine the most critical lines in the network. In particular, $N-1$ contingency is applied to the test cases to find the loss of which line yields to higher LODF values.  The $N-2$ contingency is also investigated in this paper. Because of the computational complexity of $N-2$ contingency analysis for larger test cases, this study only considers smaller test cases for analysis. Graphlet analysis of those lines is performed in order to find a relationship between lines with high LODFs and the local structure of network. It is explored whether any critical lines contain certain sorts of graphlets. This will expose any vulnerable network components that have the potential to cause a cascading failure. That will eventually provide power system operators with crucial recommendations to improve power system vulnerability in the event of a cascading failure.
The following is how this document is structured. Section~\ref{Overview of Line Outage Distribution Factor} reviews LODF metric. Section~\ref{Power Grid as a Complex Network} explains how power networks can be interpreted as complex network. Section~\ref{Methodology} presents our methodology. Section~\ref{Results} shows the results and Section~\ref{Conclusion} concludes the paper.

\section{Overview of Line Outage Distribution Factor}
\label{Overview of Line Outage Distribution Factor}

This section explains the concept of the line outage distribution factors (LODFs) that is employed in this study. To compute LODFs, we first need to consider Injection Shift Factors (ISF). The injection shift factor is a linear sensitivity factor which measures the sensitivities of active power line flows due to injection of active power bus injections~\cite{chen2016generalized}. Assume a power system consists of $N$ buses and $L$ lines. The injection shift factor $\iota^k_{lm}$ of line $lm$ represents the sensitivity of the change in power flow in line $lm$ due to a change in real power injection at bus $k \in N$. The injection shift factor matrix of a power system can be computed by DC power flow.

\begin{align}
\label{isf matrix}
\underline{\iota} = {\underline{B}^{-1}_n \underline{A} \underline{B_l} }
\end{align}

Here, $\underline{B}^{-1}_n  \in \mathbb{R}^{L\times L}$ represents the original branch susceptance matrix,  $\underline{A}  \in \mathbb{R}^{L\times N}$ denotes reduced incidence matrix, and $\underline{B_l}  \in \mathbb{R}^{L\times N}$  represents reduced nodal susceptance matrix. Power transmission distribution factor (PTDF) for line $lm$ for incident $z$ of real power change $\Delta P$ in bus $c$ to bus $d$ can be defined as,

\begin{align}
\label{ptdf}
{\tau^{z}_{lm}} = {\iota^c_{lm} - \iota^d_{lm} }
\end{align}

Consequently, change of power flow in line $lm$ for $\Delta P$ power flow change in bus $c$ to bus $d$ can be defined in terms of power as,

\begin{align}
\label{ptdf_power}
{\Delta\omega^{z}_{l_m}} = {\tau^{z}_{l_m}} \Delta P
\end{align}

PTDF matrix of the entire system can be obtained using~\eqref{isf matrix}. 
\begin{align}
\label{ptdf_matrix}
{\underline\tau_{L}} = {\underline{ \iota A^T}}
\end{align}
For outage of line $lo$, LODF of line $lm$ can be considered as ${\phi}^{lo}_{lm}$. If pre-outage real power flow in line $lo$ is $\omega_{lo}$, change of power flow ${\Delta\omega^{z(lo)}_{lm}}$ in line $lm$ after the outage of line $lo$ can be written in terms of LODF as follows~\cite{guler2007generalized}

\begin{align}
\label{lodf}
{\Delta\omega^{z(lo)}_{lm} = {\phi}^{lo}_{lm}{\omega_{lo}}}
\end{align}

Using equation~\eqref{ptdf_power}, equation~\eqref{lodf} can be extended as.
\begin{align}
\label{lodf_ext}
{\phi}^{lo}_{lm} = {\frac{{\tau^{z(lo)}_{lm}} }{1-\tau^{z(lo)}_{lo}}}
\end{align}
In this paper, we first determine the most critical lines using LODFs. For a single test case consisting of $N$ buses and $L$ lines, the LODF matrix will be of size $L \times L$. In particular, for a single line outage, there will be $L-1$ number of LODFs. The maximum LODF value is recorded among those $L-1$ LODFs. We conduct the $N-1$ contingency analysis for all lines in the network to determine maximum LODF for all contingencies. After finding LODF matrix, a limited number of lines with highest and lowest LODF are investigated from local structure point of view. A similar approach is taken for $N-2$ contingency. Following the outage of a line, it is determined that the outage of whichever line follows, in conjunction with the first, stresses the network the most. As a result, some combinations of lines for $N-2$ contingency can be obtained. Each combination gives a maximum LODF for an $N-2$ contingency. In this process, LODF helps to determine the most critical lines of a network in both $N-1$ and $N-2$ contingencies. In following sections, the local network characteristics of those critical lines are investigated.

\section{Power Grid as a Complex Network}
\label{Power Grid as a Complex Network}

Network Graphlets are building blocks of networks. Analysis of Graphlets is found to be an indispensable tool for understanding local network structure, in contrast to measures based on node degree distribution and its functions that primarily address a global network topology~\cite{ dey2017motif }. Examining the local structure of any grid requires power grid to be modeled as a complex network. The grid can be represented as an undirected graph $G(V,E)$, where $V$ represents the grid's buses and $E$ represents the grid's lines. A graph $G'(V',E')$ is a subgraph of graph ${G}$, ${G'} \subseteq {G}$, only when ${V'} \subseteq {V}$ and $ {E'} \subseteq {E}$.  The subgraph $G'$ can be called an induced subgraph of $G$ if $E'$ contains all edges $e_{uv}\in E$ such that $u,v\in V'$. Graphs $G'$ and $G''$ can be called isomorphic if there exists a bijection $h:V'\rightarrow V''$ such that any two adjacent nodes $u,v \in V' $ of $G'$ are also adjacent in $G''$ after the mapping occurs. If $G_k = (V_k,E_k)$ is a k-node subgraph of G and if there exists an isomorphism between $G_k$ and $G'$ where $G' \in G$, then there is an occurrence of $G_k$ in $G$. In this paper, we only  investigate 4-node subgraphs. To this end, we use the FANMOD algorithm, represented in Algorithm~\ref{AlgorithmI}, to enumerate subgraphs~\cite{wernicke2006fanmod, wernicke2006efficient}. In this algorithm, first all 4-node subgraphs are enumerated and then classified. Here, the notation of $N(\{v\})$ means neighbours of nodes $v \in V$. $V_{subgraph}$ means set of four node connected vertices. Furthermore, if vertex $w \in V/ V_{subgraph}$ , then $N_{excl}(w,V_{subgraph})$ denotes the exclusive neighbourhood of vertex $w$, in which no vertices belong to either $V_{subgraph}$ or the neighbourhood of $V_{subgraph}$. After enumerating subgraphs, we apply subgraph isomorphism to classify all 4-node subgraphs.

\begin{algorithm}
\caption{$Enumerate Subgraph$}
\label{AlgorithmI}
\begin{algorithmic}
\STATE \textbf{Input}: Graph $G(V,E)$
\STATE \textbf{Output}: All 4-nodes subgraph in graph $G$
\FOR{each vertex  $v \in V$}
\STATE $V_{extension} \leftarrow {u \in N({\{ v \}}): u>v }$
\STATE \textbf{Call} $ExtendSubgraph$(${\{v\}},V_{extension},v$)
\ENDFOR
\end{algorithmic}
\begin{algorithmic}
\STATE $ExtendSubgraph$(${\{v\}},V_{extension},v$)
\IF{$|V_{subgraph}|=4$}
\STATE \textbf{Output: $G[V_{subgraph}]$} 
\ENDIF   
\WHILE{$V_{extension}\neq \emptyset$}
\STATE Remove an arbitrary chosen vertex w from $V_{extension}$
\STATE $V'_{\text{extension}} \leftarrow V_{\text{extension}} \cup \{ u \in N_{\text{excl}}(w,V_{\text{subgraph}}) : u > v \}$
\STATE\textbf{Call}  ExtendSubgraph$(V_{subgraph} \cup {\{w\}},V'_{extension},v$)
\ENDWHILE
\end{algorithmic}
\end{algorithm} 

\begin{figure}
\centering
\captionsetup{justification=centering,margin=6cm}
{\includegraphics[width=5.4cm]{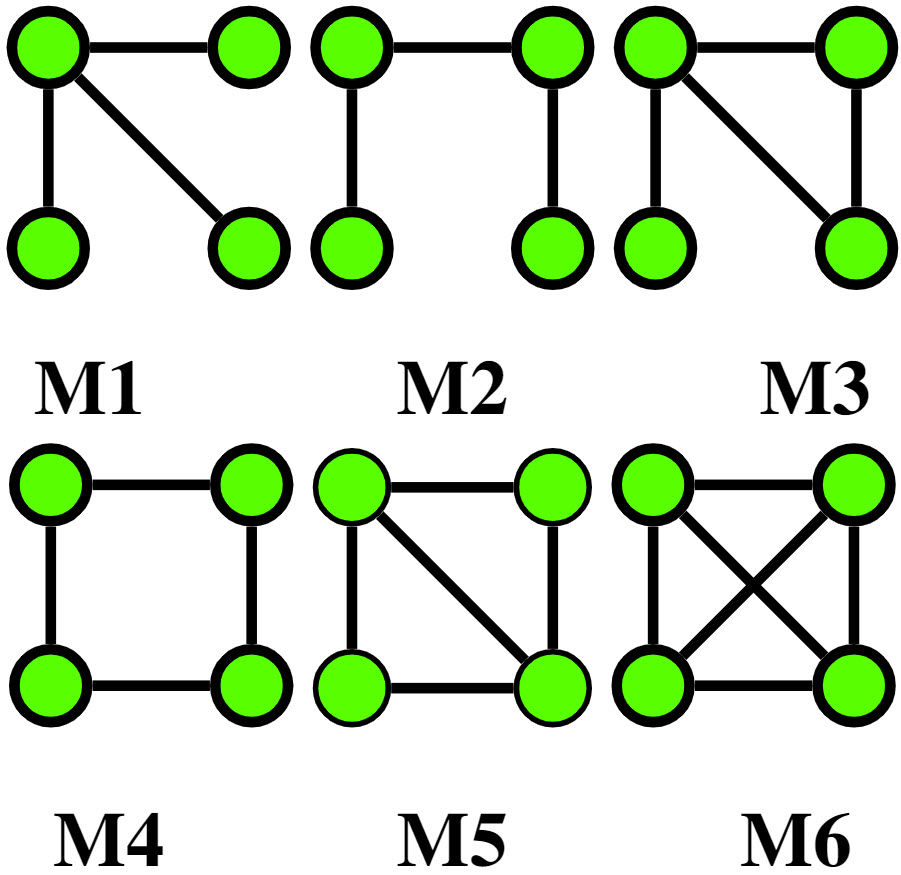}}
\caption{Different type of four node graphlets. In a power grid network it is very usual to find subgraphs M1 and M2}
\label{fig:motif_types}
\end{figure}

In general, there are six types of four-node connected subgraphs that can be found in any power grid network. Fig.~\ref{fig:motif_types} shows four-node subgraphs which are labeled as M1-M6. Subgraphs M1, M2 and M3 represent radial structures, while M4, M5 and M6 symbolize ring and mesh structures. Different types of subgraph analysis have been done in the literature. 
There have been much research on graphlet analysis in a network. However, the impact of higher-order structures such as graphlets in power grids' network is less studied.  In this study, our methodology focuses on how many times a graphlet can be found on a particular line.

\section{Methodology}
\label{Methodology}
\subsection{N-1 Contingency}
\label{N-1 Contin}
\begin{figure}
\centering
\captionsetup{justification=centering,margin=6cm}
{\includegraphics[width=6.4cm]{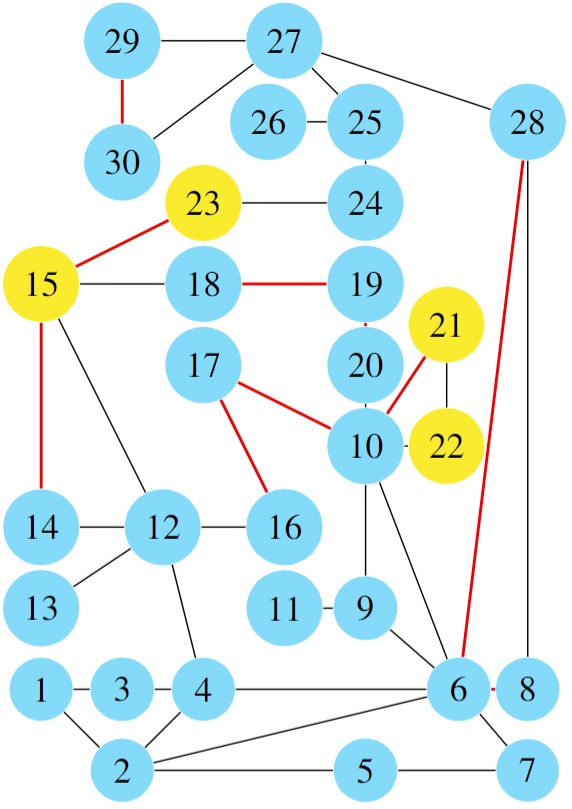}}
\caption{IEEE 30-bus system with red lines marked as the most critical lines based on $N-1$ contingency LODF analysis. Lines 15-23 and 21-22 represents $N-2$ contingency combination of lines that produce highest LODF impact.}
\label{fig:N-1 contingency}
\end{figure}

The LODF calculates the redistribution of a line's load among the other lines. However, we observe that for outages of a few particular lines, the LODF of some other lines are larger. This essentially means that outages on those critical lines are putting a strain on the network. In other words, for the outage of those critical lines, LODF impacts on other lines are significantly higher. In this section, the graphlet characteristics of those critical lines is observed. In particular, we observe how many times a particular graphlet can be found on a line. For $N-1$ contingency, the most critical lines of the IEEE 30-bus network based on LODF impact are illustrated by red lines in Fig 2. Based on the LODF analysis, it can be said that outage on those lines strain the network most. From a graph theory point of view, graphlet analysis can be used to determine which 4-node graphlets those critical lines are belong to. As a consequence, the local structural properties of those critical lines will be revealed. In this connection, we determine how many times a particular graphlet is incident on a line. Table~\ref{table:Subgraph characteristics of most critical lines} and Table \ref{table:Subgraph characteristics of least critical lines} represent graphlet characteristics of most and least critical lines based on $N-1$ LODF analysis for IEEE 30-bus network. When two tables are compared, it is clear that the local structural characteristics of the two distinct type of lines differ noticeably. In the table of most critical lines, there are no M3-M6 graphlets in any lines. On the contrary, there is an M4 ring graphlet present for all lines in the table of least critical lines. In other words, if a line is part of a ring graphlet M4, its outage is unlikely to have a significant influence on the network. In addition to that, from figure \ref{fig:N-1 contingency}, it is also clear that there are no critical lines in the graph region where the M4 ring graphlet exists. As a result, the area is considerably stabilized by the presence of a ring structure. On the contrary, most of the critical lines are part of some radial structured graphlet such as M1 and M2. In the following section, we will observe the local structural properties of lines with the highest LODF impact among different networks.

\begin{table}
    \centering
    \caption{graphlet characteristics of most critical lines in IEEE 30-bus test system}
    \begin{tabular}{|l|l|l|l|l|l|l|l|}
    \hline
        From & To & M1 & M2 & M3 & M4 & M5 & M6 \\ \hline
        14 & 15 & 1 & 2 & 5 & 0 & 0 & 0 \\ \hline
        15 & 23 & 2 & 9 & 1 & 0 & 0 & 0 \\ \hline
        16 & 17 & 0 & 10 & 0 & 0 & 0 & 0 \\ \hline
        29 & 30 & 0 & 0 & 2 & 0 & 0 & 0 \\ \hline
        10 & 21 & 5 & 8 & 6 & 0 & 0 & 0 \\ \hline
        10 & 17 & 3 & 2 & 1 & 0 & 0 & 0 \\ \hline
        19 & 20 & 0 & 7 & 0 & 0 & 0 & 0 \\ \hline
        18 & 19 & 0 & 5 & 0 & 0 & 0 & 0 \\ \hline
        6 & 28 & 4 & 14 & 7 & 0 & 0 & 0 \\ \hline
        6 & 8 & 4 & 8 & 7 & 0 & 0 & 0 \\ \hline
    \end{tabular}
    \label{table:Subgraph characteristics of most critical lines}
\end{table}

 \begin{table}
    \centering
    \caption{graphlet characteristics of least critical lines in IEEE 30-bus test system}
    \begin{tabular}{|l|l|l|l|l|l|l|l|}
    \hline
        From & To & M1 & M2 & M3 & M4 & M5 & M6 \\ \hline
        5 & 7 & 0 & 7 & 0 & 1 & 0 & 0 \\ \hline
        2 & 4 & 2 & 8 & 9 & 1 & 0 & 0 \\ \hline
        6 & 7 & 12 & 14 & 3 & 1 & 0 & 0 \\ \hline
        2 & 5 & 2 & 9 & 1 & 1 & 0 & 0 \\ \hline
        2 & 6 & 9 & 16 & 11 & 1 & 0 & 0 \\ \hline
        3 & 4 & 2 & 12 & 1 & 1 & 0 & 0 \\ \hline
        1 & 3 & 0 & 4 & 0 & 1 & 0 & 0 \\ \hline
        1 & 2 & 2 & 9 & 1 & 1 & 0 & 0 \\ \hline
    \end{tabular}
   \label{table:Subgraph characteristics of least critical lines}
\end{table}

\subsection{N-2 Contingency}
\label{N-2 Contingency}

In $N-2$ contingency, we observe network performance in the case of an outage of two lines. In $N-1$ contingency for a network of $L$ lines,  we have a LODF matrix of  $L \times \ L$ size. However, for $N-2$ contingency, we have $L$ number of $L-1$ sized matrices. From each of those $L-1$ sized matrices, we choose the row with highest LODF value. Consequently, we get $L$ combination of two lines which actually represents lines from $N-2$ contingency. In particular, that means we first simulate the outage of a single line, then we observe the outage of which other line strains the network most. In that way, we get $L$ combination of $N-2$ contingency lines and then we determine which single combination yields the maximum LODF impact. In Fig~\ref{fig:N-1 contingency}, line combinations $15-23$ and $21-22$ has the highest LODF among other combinations. If we look at the network from a structural point of view, we observe that those two lines do not belong to any ring or meshed graphlets. In the next section, we present the results of highest LODF impact line combination from $N-2$ contingency among different networks.

\section{Results}
\label{Results}

In the previous section, we observed graphlet characteristics of the IEEE 30-bus network. Here, we apply the proposed approach on ``pglib$\_$793$\_$goc'' and ``pglib$\_$1354$\_$pegase'' as large networks from the PGLib-OPF benchmark library~\cite{babaeinejadsarookolaee2019power}. For ``pglib$\_$793$\_$goc'' test case, we observed that it has no M3, M4, M5, M6 graphlets in its most critical lines, while in it's least critical lines, M4 graphlets can be found. Similarly, for the case ``pglib$\_$1354$\_$pegase'', we noticed that, there is a very small number of M3, M4, M5, M6 graphlets in it's most critical lines. In contrast, most of the least critical lines in this test case are belong to M3 and M4 graphlets. 

\begin{table*}
    \centering
    \caption{Graphlet characteristics of Most and least critical lines of ``pglib$\_$793$\_$goc'' test system}
    \begin{tabular}{|l|l|l|l|l|l|l|l|l|l|l|l|}
    \hline
        \multicolumn{5}{c}{Most Critical Edges} && \multicolumn{5}{c}{Least Critical Edges} \\ \hline
        M1 & M2 & M3 & M4 & M5 & M6 & M1 & M2 & M3 & M4 & M5 & M6 \\ \hline
        1 & 9 & 0 & 0 & 0 & 0 & 0 & 2 & 0 & 1 & 0 & 0 \\ \hline
        4 & 11 & 0 & 0 & 0 & 0 & 3 & 11 & 0 & 1 & 0 & 0 \\ \hline
        3 & 11 & 0 & 0 & 0 & 0 & 6 & 14 & 0 & 1 & 0 & 0 \\ \hline
        1 & 3 & 0 & 0 & 0 & 0 & 3 & 6 & 0 & 1 & 0 & 0 \\ \hline
        1 & 7 & 0 & 0 & 0 & 0 & 3 & 8 & 0 & 1 & 0 & 0 \\ \hline
        2 & 15 & 0 & 0 & 0 & 0 & 0 & 4 & 0 & 1 & 0 & 0 \\ \hline
        7 & 19 & 0 & 0 & 0 & 0 & 2 & 9 & 0 & 1 & 0 & 0 \\ \hline
        1 & 6 & 0 & 0 & 0 & 0 & 12 & 29 & 0 & 1 & 0 & 0 \\ \hline
        6 & 12 & 0 & 0 & 0 & 0 & 2 & 12 & 0 & 2 & 0 & 0 \\ \hline
        1 & 11 & 0 & 0 & 0 & 0 & 1 & 5 & 0 & 1 & 0 & 0 \\ \hline
    \end{tabular}
\end{table*}

\begin{table*}
    \centering
    \caption{Graphlet characteristics of Most and least critical lines of ``pglib$\_$1354$\_$pegase'' test system}
    \begin{tabular}{|l|l|l|l|l|l|l|l|l|l|l|l|}
    \hline
        \multicolumn{5}{c}{Most Critical Edges} && \multicolumn{5}{c}{Least Critical Edges} \\ \hline
        M1 & M2 & M3 & M4 & M5 & M6 & M1 & M2 & M3 & M4 & M5 & M6 \\ \hline
        36 & 9 & 1 & 0 & 0 & 0 & 3 & 13 & 0 & 1 & 0 & 0 \\ \hline
        1 & 3 & 0 & 0 & 0 & 0 & 3 & 11 & 0 & 1 & 0 & 0 \\ \hline
        3 & 6 & 0 & 0 & 0 & 0 & 9 & 16 & 1 & 2 & 0 & 0 \\ \hline
        3 & 6 & 0 & 0 & 0 & 0 & 10 & 20 & 1 & 2 & 0 & 0 \\ \hline
        6 & 22 & 0 & 0 & 0 & 0 & 10 & 20 & 1 & 2 & 0 & 0 \\ \hline
        13 & 38 & 0 & 0 & 0 & 0 & 26 & 29 & 2 & 1 & 0 & 0 \\ \hline
        5 & 14 & 0 & 0 & 0 & 0 & 3 & 9 & 1 & 1 & 0 & 0 \\ \hline
        10 & 18 & 0 & 0 & 0 & 0 & 52 & 28 & 16 & 1 & 0 & 0 \\ \hline
        11 & 27 & 0 & 0 & 0 & 0 & 4 & 27 & 0 & 2 & 0 & 0 \\ \hline
        2 & 15 & 0 & 0 & 0 & 0 & 7 & 26 & 0 & 1 & 0 & 0 \\ \hline
    \end{tabular}
\end{table*}

\begin{figure}
\centering
\captionsetup{justification=centering,margin=6cm}
{\includegraphics[width= 7 cm, trim={0.0cm 0.0cm 0.0cm 0.0cm},clip]{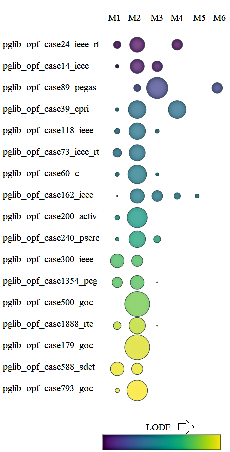}}
\caption{Relation between graphlet characteristics of different networks and their maximum $N-1$ LODF impacts. Size of the bubble represents graphlet percentages and color of the bubbles changes with Maximum $N-1$ LODF impact possible in each cases}
\label{fig: N-1 motifs}
\end{figure}

In order to generalize the finding, we investigate the proposed approach on different test cases with various sizes. To this end, we first determined the most critical lines of different test cases using $N-1$ LODF analysis. We then determined which graphlets they feature and the quantity of each graphlet present. In particular, total graphlet count of the most important line is determined. Then it is measured which graphlets account for what percentage of the entire graphlet count. This gives us the graphlet characteristics of that line. By taking the percentage of graphlets in a line, we can compare the results across different networks. The maximum LODF impact of most critical line of a network varies from network to network. Our objective is to relate the maximum LODF for a network possible and subsequent graphlet characteristics of the line for the outage of which the network experiences the maximum LODF. We present our findings in Fig.~\ref{fig: N-1 motifs}, where the size of the bubble represents each graphlet percentage of the most important line of a network. Also, the color of the bubbles changes with the  maximum LODF of that most important line. From Fig.~\ref{fig: N-1 motifs}, it is evident that, cases that contain very high LODF impact line, contain no or very small percentages of ring or meshed graphlets such as M3, M4, M5 and M6 in their most critical line. Conversely, there are relatively higher percentages of ring or meshed graphlets in cases that have relatively low LODF impact lines.

\begin{figure}
\centering
\captionsetup{justification=centering,margin=6cm}
{\includegraphics[width= 7 cm,trim={0.0cm 0.2cm 0.4cm 0.2cm},clip]{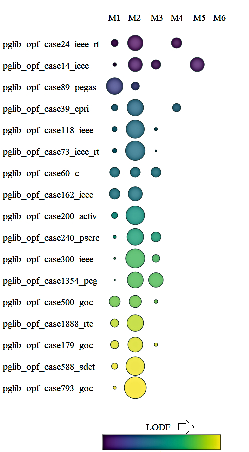}}
\caption{Relation between graphlet characteristics of different networks and their maximum $N-2$ LODF impacts. Size of the bubble represents graphlet percentages and color of the bubbles changes with Maximum $N-2$ LODF impact possible in each cases}
\label{fig: N-2 motifs}
\end{figure}

From Section \ref{Methodology}, we get $L$ combination of two $N-2$ contingency lines for a network consisted of $L$ lines. The first line corresponds to $N-1$ contingency and the second line corresponds to $N-2$ contingency. For most combinations, we can observe certain $N-2$ lines recurring. However, we focus on the $N-2$ contingency line for the combination that produces the maximum LODF impact. Fig.~\ref{fig: N-2 motifs} shows the graphlet analysis of lines for $N-2$ contingency. From Fig.~\ref{fig: N-2 motifs} it can be noticed that, cases with very high LODF impact lines have less percentage of ring or meshed graphlets such as M3, M4, M5 and M6 in their most critical lines. This results can be easily extended to higher order contingency analysis and help power system operator to focus on specific part of the network for finding critical lines. Finding graphlets for larger test cases can be done offline which helps system operators to quickly determine potential critical lines in the system. 

\section{Conclusion}
\label{Conclusion}
Graph theoretical analysis has been extensively used to investigate power grid vulnerabilities. However, graphlet analysis in tandem with grid vulnerabilities is less studied. To address this gap, this paper first determines the most important line for $N-1$ contingency based on the LODF analysis. Then, it focuses on the graphlet characteristics of that critical line. This strategy is applied to various test cases and graphlet characteristics of the most important lines of those networks are investigated. This provides a comparison of the highest LODF possible for a network and graphlet features of the line for outages in which that LODF occurs. It is revealed that networks with a higher percentage of ring or meshed graphlets such as M3, M4, M5, M6 on their most important line experience less LODF for outage of that line. On the contrary, cases with very high LODF for outage on its most important line have very little to no percentages of ring or meshed graphlets in that line. This study also expands to, $N-2$ contingency where the same type of scenario is observed, which actually emphasizes the importance of ring or meshed graphlets in power grid networks for more resilient grids. The proposed graphlet analysis helps the power system operator to quickly determine possible candidates for critical lines.   
\bibliographystyle{IEEEtran}
\bibliography{ref}

\begin{thebibliography}{10}
\providecommand{\url}[1]{#1}
\csname url@samestyle\endcsname
\providecommand{\newblock}{\relax}
\providecommand{\bibinfo}[2]{#2}
\providecommand{\BIBentrySTDinterwordspacing}{\spaceskip=0pt\relax}
\providecommand{\BIBentryALTinterwordstretchfactor}{4}
\providecommand{\BIBentryALTinterwordspacing}{\spaceskip=\fontdimen2\font plus
\BIBentryALTinterwordstretchfactor\fontdimen3\font minus
  \fontdimen4\font\relax}
\providecommand{\BIBforeignlanguage}[2]{{%
\expandafter\ifx\csname l@#1\endcsname\relax
\typeout{** WARNING: IEEEtran.bst: No hyphenation pattern has been}%
\typeout{** loaded for the language `#1'. Using the pattern for}%
\typeout{** the default language instead.}%
\else
\language=\csname l@#1\endcsname
\fi
#2}}
\providecommand{\BIBdecl}{\relax}
\BIBdecl

\bibitem{galiana1984bound}
F.~Galiana, ``Bound estimates of the severity of line outages in power system
  contingency analysis and ranking,'' \emph{IEEE Transactions on Power
  Apparatus and Systems}, no.~9, pp. 2612--2624, 1984.

\bibitem{liang20162015}
G.~Liang, S.~R. Weller, J.~Zhao, F.~Luo, and Z.~Y. Dong, ``The 2015 ukraine
  blackout: Implications for false data injection attacks,'' \emph{IEEE
  Transactions on Power Systems}, vol.~32, no.~4, pp. 3317--3318, 2016.

\bibitem{Boyaci_propagation}
O.~Boyaci, M.~R. Narimani, K.~Davis, and E.~Serpedin, ``Spatio-temporal failure
  propagation in cyber-physical power systems,'' in \emph{2022 3rd
  International Conference on Smart Grid and Renewable Energy (SGRE)}, 2022,
  pp. 1--6.

\bibitem{chen2007identification}
X.~Chen, K.~Sun, Y.~Cao, and S.~Wang, ``Identification of vulnerable lines in
  power grid based on complex network theory,'' in \emph{2007 IEEE power
  engineering society general meeting}.\hskip 1em plus 0.5em minus 0.4em\relax
  IEEE, 2007, pp. 1--6.

\bibitem{coelho2018new}
E.~P.~R. Coelho, M.~H.~M. Paiva, M.~E.~V. Segatto, and G.~Caporossi, ``A new
  approach for contingency analysis based on centrality measures,'' \emph{IEEE
  Systems Journal}, vol.~13, no.~2, pp. 1915--1923, 2018.

\bibitem{gao2019identification}
Q.~Gao, Y.~Wang, X.~Cheng, J.~Yu, X.~Chen, and T.~Jing, ``Identification of
  vulnerable lines in smart grid systems based on affinity propagation
  clustering,'' \emph{IEEE Internet of Things Journal}, vol.~6, no.~3, pp.
  5163--5171, 2019.

\bibitem{tejada2017security}
D.~A. Tejada-Arango, P.~S{\'a}nchez-Mart{\i}n, and A.~Ramos, ``Security
  constrained unit commitment using line outage distribution factors,''
  \emph{IEEE Transactions on power systems}, vol.~33, no.~1, pp. 329--337,
  2017.

\bibitem{narimani2020tightening}
M.~R. Narimani, D.~K. Molzahn, and M.~L. Crow, ``Tightening qc relaxations of
  ac optimal power flow problems via complex per unit normalization,''
  \emph{IEEE Transactions on Power Systems}, vol.~36, no.~1, pp. 281--291,
  2020.

\bibitem{narimani2018comparison}
M.~R. Narimani, D.~K. Molzahn, H.~Nagarajan, and M.~L. Crow, ``Comparison of
  various trilinear monomial envelopes for convex relaxations of optimal power
  flow problems,'' in \emph{2018 IEEE Global Conference on Signal and
  Information Processing (GlobalSIP)}.\hskip 1em plus 0.5em minus 0.4em\relax
  IEEE, 2018, pp. 865--869.

\bibitem{narimani2018improving}
M.~R. Narimani, D.~K. Molzahn, and M.~L. Crow, ``Improving qc relaxations of
  opf problems via voltage magnitude difference constraints and envelopes for
  trilinear monomials,'' in \emph{2018 Power Systems Computation Conference
  (PSCC)}.\hskip 1em plus 0.5em minus 0.4em\relax IEEE, 2018, pp. 1--7.

\bibitem{narimani2018empirical}
M.~R. Narimani, D.~K. Molzahn, D.~Wu, and M.~L. Crow, ``Empirical investigation
  of non-convexities in optimal power flow problems,'' in \emph{2018 Annual
  American Control Conference (ACC)}.\hskip 1em plus 0.5em minus 0.4em\relax
  IEEE, 2018, pp. 3847--3854.

\bibitem{narimani2020strengthening}
M.~R. Narimani, \emph{Strengthening QC relaxations of optimal power flow
  problems by exploiting various coordinate changes}.\hskip 1em plus 0.5em
  minus 0.4em\relax Missouri University of Science and Technology, 2020.

\bibitem{narimani2021generalized}
M.~R. Narimani, H.~Huang, A.~Umunnakwe, Z.~Mao, A.~Sahu, S.~Zonouz, and
  K.~Davis, ``Generalized contingency analysis based on graph theory and line
  outage distribution factor,'' \emph{IEEE Systems Journal}, vol.~16, no.~1,
  pp. 626--636, 2021.

\bibitem{huang2021toward}
H.~Huang, Z.~Mao, M.~R. Narimani, and K.~R. Davis, ``Toward efficient wide-area
  identification of multiple element contingencies in power systems,'' in
  \emph{2021 IEEE Power \& Energy Society Innovative Smart Grid Technologies
  Conference (ISGT)}.\hskip 1em plus 0.5em minus 0.4em\relax IEEE, 2021, pp.
  01--05.

\bibitem{milo2002network}
R.~Milo, S.~Shen-Orr, S.~Itzkovitz, N.~Kashtan, D.~Chklovskii, and U.~Alon,
  ``Network motifs: simple building blocks of complex networks,''
  \emph{Science}, vol. 298, no. 5594, pp. 824--827, 2002.

\bibitem{dey2017motif}
A.~K. Dey, Y.~R. Gel, and H.~V. Poor, ``Motif-based analysis of power grid
  robustness under attacks,'' in \emph{2017 IEEE Global Conference on Signal
  and Information Processing (GlobalSIP)}.\hskip 1em plus 0.5em minus
  0.4em\relax IEEE, 2017, pp. 1015--1019.

\bibitem{chen2016generalized}
Y.~C. Chen, S.~V. Dhople, A.~D. Dom{\'\i}nguez-Garc{\'\i}a, and P.~W. Sauer,
  ``Generalized injection shift factors,'' \emph{IEEE Transactions on Smart
  Grid}, vol.~8, no.~5, pp. 2071--2080, 2016.

\bibitem{guler2007generalized}
T.~Guler, G.~Gross, and M.~Liu, ``Generalized line outage distribution
  factors,'' \emph{IEEE Transactions on Power systems}, vol.~22, no.~2, pp.
  879--881, 2007.

\bibitem{wernicke2006fanmod}
S.~Wernicke and F.~Rasche, ``Fanmod: a tool for fast network motif detection,''
  \emph{Bioinformatics}, vol.~22, no.~9, pp. 1152--1153, 2006.

\bibitem{wernicke2006efficient}
S.~Wernicke, ``Efficient detection of network motifs,'' \emph{IEEE/ACM
  transactions on computational biology and bioinformatics}, vol.~3, no.~4, pp.
  347--359, 2006.

\bibitem{babaeinejadsarookolaee2019power}
S.~Babaeinejadsarookolaee, A.~Birchfield, R.~D. Christie, C.~Coffrin,
  C.~DeMarco, R.~Diao, M.~Ferris, S.~Fliscounakis, S.~Greene, R.~Huang
  \emph{et~al.}, ``The power grid library for benchmarking ac optimal power
  flow algorithms,'' \emph{arXiv preprint arXiv:1908.02788}, 2019.

\end{thebibliography}
\end{document}